\documentclass[a4paper,12pt]{article}
\usepackage{mathtext}
\usepackage{graphicx}
\usepackage[T2A]{fontenc}
\usepackage[utf8]{inputenc}
\usepackage[english]{babel}
\usepackage{latexsym,amsfonts,amssymb,amsmath,longtable,amsthm}
\usepackage{cite}
\usepackage[pic]{xy}
\usepackage{bbm,makeidx}
\usepackage[centerlast,small]{caption}
\usepackage{amsmath}
\begin{document}

\begin{center}
{\LARGE Twisted conjugacy classes in Chevalley groups\footnote{The author was supported by Russian foundation of basic research (project code 12-01-33102).
}}

~

~

{\large Timur R. Nasybullov}

~

Sobolev Institute of Mathematics \& Novosibirsk State University

ntr@math.nsc.ru
\end{center}

\newcounter{thelem}
\newtheorem{lem}{{\scshape Lemma}}
\newtheorem{ttt}{{\scshape Theorem}}
\newtheorem{exa}{{\scshape Example}}
\newtheorem{cor}{{\scshape Corollary}}
\newtheorem{prp}[lem]{{\scshape Proposition}}

~~~~~~~~~~\parbox{11cm}{\begin{center}\textbf{Abstract}
\end{center}
\small We prove that Chevalley group over the field $F$ of zero characteristic possess $R_{\infty}$ property, if $F$ has torsion group of automorphisms or $F$ is an algebraically closed field which has finite transcendence degree over $\mathbb{Q}$.  As a consequence we obtain that the twisted conjugacy class $[e]_{\varphi}$ of unit element is a subgroup of Chevalley group if and only if $\varphi$ is central automorphism.\normalsize}

\section{Introduction}

Let $G$ be a group and $\varphi : G \longrightarrow G$ be an automorphism of $G$. Elements $x$ and $y$ of group $G$ are said to be $ \emph{twisted} ~\varphi$-\emph{conjugated} or simply $\varphi$-\emph{conjugated} ($x\sim_{\varphi}y$), if there exists an element $z$ of group $G$, such that the equality $x=zy\varphi(z^{-1})$ holds. If $\varphi$ is an identical automorphism, than we have the definition of conjugated elements here. The relation of $\varphi$-conjugation is an equivalence relation and here we can speak about $\varphi$-conjugacy classes. We use the symbol $[x]_{\varphi}$ to denote $\varphi$-conjugacy class of the element $x$. The number $R(\varphi)$ of $\varphi$-conjugacy classes is called \emph{Reidemeister number} of $\varphi$.

The classical Bernside theorem \cite[\S 10, theorem 2]{kir} states that the number of conjugacy classes of finite group $G$ is equal to the number of its irreducible complex (and therefore unitary) representations. Currently analogue of this theorem, which is called Twisted Burnside-Frobenius Theorem (TBFT) actively studied. TBFT has its origin in
the following cojecture of A.~Fel'shtyn and R.~Hill \cite{FH}: the number $R(\varphi)$ is equal to the number of fixed points of the map $\widehat{\varphi}$, which is induced by $\varphi$ onto the set of equivalence classes of unitary irreducible representations of $G$. The first step for solving this problem is to describe such a groups, that the Reidemeister number  $R(\varphi)$ is infinite for any its automorphism $\varphi$. We say that such groups possess $R_{\infty}$ property.

The problem of determining groups which possess $R_{\infty}$ property was formulated by A. Fel'shtyn and R. Hill \cite{FH}. Non-elementary Gromov hyperbolic groups \cite{F,LL,FG1}, Baumslag-Soliter groups $BS(m,n)$ for $(m,n)\neq(1,1)$ \cite{FG}, some free nilpotent and free solvable groups \cite{KR1,KR} are known to possess $R_{\infty}$ property.
In the paper \cite{Nas} this question is investigated for  general and special linear groups. It is stated, that the general linear group ${\rm GL}_n(K)$ and the special linear group ${\rm SL}_n(K)$ ($n\geq 3$) possess $R_{\infty}$ property if $K$ is an infinite integral domain with trivial automorphism group, or $K$ is an integral domain of zero characteristic with torsion automorphism group ${\rm Aut}~K$.

In present paper we study Chevalley groups over different field. This class of linear groups contains all the classical linear groups, in particular, special linear group and symplectic group. Here we prove that Chevalley group over algebraically closed field $F$ of zero characteristic possess $R_{\infty}$ property if the transcendence degree of the field $F$ over  $\mathbb{Q}$ is finite (theorem 1). Also it is proved that Chevalley group over the field $F$ of zero characteristic possess $R_{\infty}$ property if automorphism group of the field $F$ is torsion (theorem 2).

The condition that the characteristic of the field $F$ is equal to zero can't be rejected. It follows from the  R.~Steinberg's result \cite[theorem 10.1]{S}, which states that for every connected linear algebraic group over algebraically closed field of non-zero characteristic there always exists such an automorphism  $\varphi$ that $R(\varphi)=1$.

 Another interest of investigating twisted conjugacy classes is to study twisted conjugacy class of the unit element $[e]_{\varphi}$.  This class contains unit element, whence the following question arises: for which groups $G$ and its automoprhisms $\varphi$ class $[e]_{\varphi}$ is a subgroup of group $G$? In the paper \cite{BarNasNes} it is stated that for automorphism $\varphi$, which acts identically modulo center of group $G$, class $[e]_{\varphi}$ is a subgroup of $G$. In present paper we prove that for Chevalley groups over some fields this result is a criterion, i. e. if $G$ is a Chevalley group over the field $F$ which is determined in the theorem 1 or in the theorem 2, than the  $\varphi$-conjugacy class of the unit element $[e]_{\varphi}$ is a subgroup of $G$ if and only if $\varphi$ acts identically modulo center of $G$ (theorems \ref{cor2} and \ref{cor3}).

Author is eternally greatfull to E. Vdovin, Ya. Nuzhhin and E. Bunina for the useful discussions and helpful notes.

\section{Preliminaries}

In this section we remind basic definitions and formulate well known results about Chevalley groups. Also we prove here some results, which are used in proofs of the main results.

The following proposition is proved in \cite[corollary 3.2]{feintro}
\begin{prp}\label{pr1} Let $\varphi$ be an automorphism of group $G$ and $\varphi_g$ be an inner automorphism of group $G$. Then $R(\varphi\varphi_g)=R(\varphi)$.
\end{prp}

The following proposition can be found in the paper \cite[lemma 2.1]{MS}.
\begin{prp}\label{pr3} Let
$$1\rightarrow N\rightarrow G \rightarrow A\rightarrow 1$$
be a shot exact sequence. If $N$ is a characteristic subgroup of $G$ and $A$ possess $R_{\infty}$ property, then $G$ possess $R_{\infty}$ property.
\end{prp}

\subsection{Chevalley groups}

All details about root systems and its properties can be found in \cite[\S3.3, 3.4]{Car},\cite[chapter II]{ham1}. Information about semisimple Lie algebras is contained in \cite{ham1}.

We fix an irreducible root system $\Phi$ with the subsystem of simple roots  $\Delta=\{\alpha_1,\dots,\alpha_l\}$. Subsystem of positive (negative) roots from $\Phi$ is denoted by $\Phi^+~(\Phi^-)$. Simple roots are enumerated by the following way.
\begin{center}
\includegraphics[width=65mm]{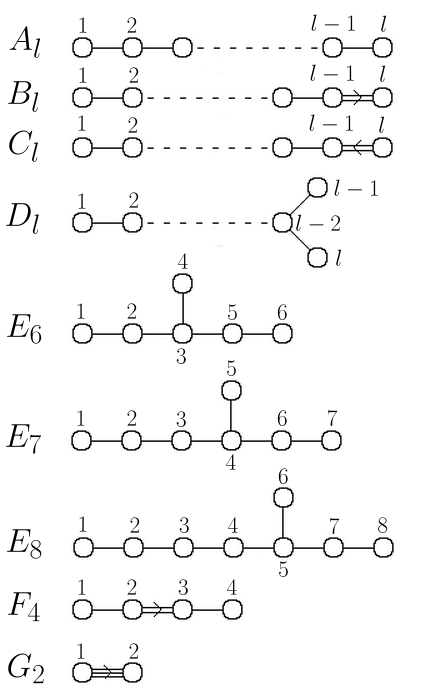}
\end{center}

Let $\mathcal{L}$ be a semisimple complex Lie algebra with Cartan subalgebra $\mathcal{H}$. Then $\mathcal{L}$ admits Cartan decomposition
$$\mathcal{L}=\mathcal{H}\oplus\sum_{\alpha\neq 0}\mathcal{L}_{\alpha},$$
where $\mathcal{L}_{\alpha}=\{x \in\mathcal{L}~|~[h,x]=\alpha(h)x~for~every~h\in \mathcal{H}\}$, and if $\mathcal{L}_{\alpha}\neq 0$, then the dimension of $\mathcal{L}_{\alpha}$ is equal to one.

 Let ${\rm ad}:\mathcal{L} \rightarrow gl(\mathcal{L})$ be an adjoint representation of  Lie algebra $\mathcal{L}$. Element ${\rm ad}x$ from $ gl(\mathcal{L})$ acts on the element $y$ from $\mathcal{L}$ by the the way ${\rm ad}x(y)=[x,y]$, where $[~,~]$ means multiplication in $\mathcal{L}$. We have bilinear form $(x,y)={\rm tr}({\rm ad}x{\rm ad}y)$ on $\mathcal{L}$, which is non-singular on $\mathcal{H}$. Therefore spaces $\mathcal{H}$ and $\mathcal{H}^*$ can be identified. All of such an elements $\alpha\in \mathcal{H}$, that $\mathcal{L}_{\alpha}\neq 0$ form the root system $\Phi$. The root system $\Phi$ and semisimple Lie algebra $\mathcal{L}$ are in one to one correspondence.

We can choose such a basis $\{h_1,\dots,h_l\}$ of $\mathcal{H}$ and such an elements $x_{\alpha}\in \mathcal{L}_{\alpha}$ for all $\alpha\in \Phi$, that the set $\{h_i,x_{\alpha}~|~\alpha\in \Phi,~i=1,\dots,l\}$ forms a basis of $\mathcal{L}$, and commutator of any two elements of this basis is an integral linear combination of the elements of the same basis. Such a basis is called \emph{Chevalley basis}.

Let us introduce elementary Chevalley groups. Let $\mathcal{L}$ be a semisimple complex Lie algebra with root system $\Phi$. In space $\mathcal{L}$ in the Chevalley basis all the maps ${\rm ad}(x_{\alpha})^k/k!$  ($k\in\mathbb{N}$) are matrices with integer entries. An
integral matrix also can be considered as a matrix over some field $F$ of zero characteristic. Let us consider the following automorphism of  $\mathcal{L}$
$$x_{\alpha}(t)={\rm exp}(t{\rm ad}(x_{\alpha}))=1+t{\rm ad}(x_{\alpha})+t^2{\rm ad}(x_{\alpha})^2/2!+\dots+t^k{\rm ad}(x_{\alpha})^k/k!+\dots$$
Since all the matrices ${\rm ad}(x_{\alpha})$ are nilpotent then this sum is finite. Automorphism $x_{\alpha}(t)$ is called \emph{elementary root element}. The subgroup of ${\rm Aut}(\mathcal{L})$, which is generated by all the automorphisms $x_{\alpha}(t)$, $\alpha \in \Phi, t\in F$ is called \emph{elementary Chevalley group} and is denoted by $\Phi(F)$.

Consider the elements
$n_{\alpha}(t)=x_{\alpha}(t)x_{-\alpha}(-t^{-1})x_{\alpha}(t)$,  $h_{\alpha}(t)=n_{\alpha}(t)n_{\alpha}(-1)$, $t\in F^*$, $\alpha\in\Phi$ of elementary Chevalley group.
The following formula, which is called commutator Chevalley formula, is proved in \cite[\S 5.2]{Car}
 $$[x_{\alpha}(t),x_{\beta}(u)]=\prod_{i\alpha+j\beta\in \Phi, i,j>0}x_{i\alpha +j\beta}(C_{ij\alpha\beta}(-t)^iu^j).$$  Here the constants $C_{ij\alpha\beta}$ don't depend on $t$ and $u$ and depend only on the type of root system. Also we have the following identities in the elementary Chevalley group $$x_{\alpha}(t_1)x_{\alpha}(t_2)=x_{\alpha}(t_1+t_2),~~~~~h_{\alpha}(t_1)h_{\alpha}(t_2)=h_{\alpha}(t_1t_2).$$

The group, which is generated by elements $\overline{x}_{\alpha}(t),~\alpha\in\Phi,t\in F$, which satisfy the same relations as elements $x_{\alpha}(t)$ of group $G=\Phi(F)$ is called \emph{universal Chevalley group} and is denoted by $\overline{G}$. In \cite[theorem 12.1.1]{Car} it is stated that for root system $\Phi\neq A_1$ we have $\overline{G}/\overline{Z}\cong G=\Phi(F)$, where $\overline{Z}$ is a center of group $\overline{G}$.
If $N$ is a central subgroup of the universal Chevalley group $\overline{G}$, then the group $\overline{G}/N$ is also called \emph{Chevalley group} of type $\Phi$ over the field $F$. For Chevalley groups we also have $(\overline{G}/N)/(Z(\overline{G}/N))\cong G=\Phi(F)$ for $\Phi\neq A_1$.

 Let us remind some facts about elementary Chevalley groups over the field $F$. The following result is proved in \cite[proposition 6.4.1]{Car}.
\begin{prp} \label{pr2} Elements $h_{\alpha}(t)$ acts on the Chevalley basis by the following way
$$h_{\alpha}(t)(h_{i})=h_{i},~~~i=1,\dots,l,$$
$$h_{\alpha}(t)(x_{\beta})=t^{A_{\alpha\beta}}x_{\beta},~~~\beta\in \Phi,$$
 where $A_{\alpha\beta}=2(\alpha,\beta)/(\beta,\beta)$ are Cartan numbers of the root system $\Phi$.
\end{prp}
Let $H$ be a subgroup of $\Phi(F)$, which is generated by all the elements $h_{\alpha}(t)$ for  $\alpha\in \Phi,~ t\in F^{*}$, i. e.
$$H=\langle h_{\alpha}(t)~|~\alpha\in \Phi,~ t\in F^{*}\rangle.$$

Let $P=\mathbb{Z}\Phi$ be the set of all integral linear combinations of elements from $\Phi$, then $P$ is an additive group, which is generated by roots. This group is a free abelian group of rank $l$  with the basis, which consists of simple roots $\Delta=\{\alpha_1,\dots,\alpha_l\}$. Homomorphism from additive group $P$ to the group $F^{*}$ is called $F$-character of $P$. $F$-character of $P$ is obviously defined by its values on simple roots.

Consider the following map from $\Phi$ to $F^*$
$$\beta \mapsto t^{A_{\alpha\beta}}, ~~~~~~\alpha \in \Phi, ~t\in F^*.$$
This map can be extended to $F$-character. Indeed, let $\chi_{\alpha,t}$ be such a map from $P$ to $F^*$, that
$$\chi_{\alpha,t}(a)=t^{2(\alpha,a)/(\alpha,\alpha)}.$$
Then $\chi_{\alpha,t}$ is $F$-character of $P$, which maps the root $\beta$ to $t^{A_{\alpha\beta}}$. All the  $F$-characters form a group with the following multiplication
$$\chi_1\chi_2(a)=\chi_1(a)\chi_2(a).$$

Using $F$-character $h(\chi)$  we can construct an automorphism of Lie algebra $\mathcal{L}$, which acts on the Chevalley basis by the following way:
$$h(\chi)(h_{i})=h_{i},~~~~~~~~~h(\chi)(x_{\beta})=\chi(\beta)x_{\beta}.$$
All the automorphisms $h(\chi)$ of Lie algebra $\mathcal{L}$ form a subgroup $\widehat{H}$ of group $\mathcal{L}$. This subgroup normalize elementary Chevalley group in the group ${\rm Aut}\mathcal{L}$. The map $\chi\mapsto h(\chi)$ is an isomorphism between group of $F$-characters and group $\widehat{H}$.
If $\chi=\chi_{\alpha,t}$, then $h(\chi)=h_{\alpha}(t)$. Hence $H$ is a subgroup of $\widehat{H}$. The following lemma gives more clear relation between groups $H$ and $\widehat{H}$.

\begin{lem}\label{lem20} For the root system  $\Phi$ consider the following set of indexes $I$ and the following polynomial $f(T)$
\begin{enumerate}
\item $\Phi=A_l, ~~~I=\{1,2,\dots,l-1\},~~~f=T^{l+1}$;
\item $\Phi=B_l, ~~~I=\{2,3,\dots,l\},~~~f=T^{2}$;
\item $\Phi=C_l, ~~~I=\{1,2,\dots,l-1\},~~~f=T^{2}$;
\item $\Phi=D_l, ~~~I=\{1,2,\dots,l-2\},~~~f=T^{2}$;
\item $\Phi=E_6, ~~~I=\{1,2,3,5\},~~~f=T^{3}$;

 $\Phi=E_7, ~~~I=\{1,2,3,4,6\},~~~f=T^{2}$;

 $\Phi=E_8, ~~~I=\{1,2,3,4,5,6,7,8\},~~~f=T$;
\item $\Phi=F_4, ~~~I=\{1,2,3,4\},~~~f=T$;
\item $\Phi=G_2, ~~~I=\{1,2\},~~~f=T$.
\end{enumerate}
Let $G=\Phi(F)$ be an elementary Chevalley group of type $\Phi$ over the field $F$ and let $h$ be an arbitrary element of group $\widehat{H}$. Then there exist elements $h_1 \in H,~h_2 \in \widehat{H}$, such that $h=h_1h_2$ and $h_2(x_{\alpha_i})=x_{\alpha_i}$ for $i \in I$.
 Moreover, if the equation $f(T)=a$ can be solved in the field $F$ for all $a$, then $H=\widehat{H}$.

\end{lem}
\textbf{Proof.} Since $H$ is an abelian group and $h_{\alpha}(t_1)h_{\alpha}(t_2)=h_{\alpha}(t_1t_2)$ for every $\alpha \in \Delta$, $t_1,t_2\in F$, then any element of group $H$ can be written
\begin{equation}\label{123}
h_{\alpha_1}(t_1)\dots h_{\alpha_l}(t_l)
\end{equation}
for some $t_i \in F^*$.

Any element $h(\chi)$ of group $\widehat{H}$ is obviously defined by its values on the simple roots, therefore it can be identified with ordered set of $l$ reversible elements of the field $F$ by the following rule
$$
h(\chi)=(\chi(\alpha_1),\dots,\chi({\alpha_l})).
$$
Then according to proposition \ref{pr2}, element $h_{\alpha}(t)$ in such denotation has the form
$$h_{\alpha}(t)=(t^{A_{\alpha\alpha_1}},t^{A_{\alpha\alpha_2}},\dots,t^{A_{\alpha\alpha_l}}).$$
Hence, by the equality (\ref{123}) any element of $H$ has the form
\begin{equation}\label{arb}
(t_1^{A_{11}}t_2^{A_{21}}\dots t_l^{A_{l1}},t_1^{A_{12}}t_2^{A_{22}}\dots t_l^{A_{l2}},\dots,t_1^{A_{1l}}t_2^{A_{2l}}\dots t_l^{A_{ll}}),
\end{equation}
where $A_{ij}=A_{\alpha_i\alpha_j}$ are Cartan numbers of root system $\Phi$.

Let $a_1,\dots,a_l$ be invertible elements of the field $F$. Consider the system of equations \begin{equation}\label{per}
T_1^{A_{1i}}T_2^{A_{2i}}\dots T_l^{A_{li}}=a_i,~~~~~i=1,\dots,l
\end{equation}
with variables $T_1,\dots,T_l$. If this system has a solution for any elements $a_1,\dots,a_l\in F^*$, then groups $H$ and $\widehat{H}$ are equal.

Really, consider an arbitrary element $h(\chi)$ of group $\widehat{H}$ and let $a_i=\chi(\alpha_i)$. Then, if $t_1,\dots,t_l$ is a solution  of (\ref{per}), then $h(\chi)=h_{\alpha_1}(t_1)\dots h_{\alpha_l}(t_l)$ and therefore  $h(\chi)\in H$.

Note that the system (\ref{per}) is equivalent to the system
\begin{align}
AT=a, \tag{\ref{per}$^\prime$}
\end{align} in the multiplicative group of the field $F$ with additively written algebraic operation. Here $A$ is a Cartan matrix of root system $\Phi$, $T=(T_1,\dots,T_l)^T$ is a column of variables, $a=(a_1,\dots,a_l)^T$ is a column of values.

Let $h=h(\chi)$ be an element of group $\widehat{H}$, which is given in the condition of lemma, and $a_1=\chi(\alpha_1), \dots, a_l=\chi(\alpha_l)$ are invertible elements of $F$. Now we are ready to prove the lemma
\begin{enumerate}
\item $\Phi=A_l$. We give detailed proof only in the case of root system  $A_l$. Proofs of other cases are similar, we only formulate some necessary facts for it.
From the statement \cite[appendix, proposition 9]{Ham} it follows, that using integral elementary transformation of rows and columns we can transform  Cartan matrix of root system $A_l$  to the diagonal matrix
$$diag(\underbrace{1,\dots,1}_{l-1}, l+1),$$
and therefore the system (\ref{per}) is solvable for any right part if the equation $T^{l+1}=a$ is solvable in $F$  for any element $a$. In this case $H=\widehat{H}$.

Otherwise we construct required elemets $h_1,h_2$.
According to the equality (\ref{arb}), arbitrary element of $H$ has form
$$(t_1^2t_2^{-1},t_1^{-1}t_2^2t_3^{-1}, t_2^{-1}t_3^2t_4^{-1},\dots, t_{l-2}^{-1}t_{l-1}^2t_l^{-1},t_{l-1}^{-1}t_l^2)$$
for some elements $t_1,\dots,t_l\in F^*$.
Let
$$\begin{cases}t_1=1,\\
t_i=1/({a_1^{i-1}a_2^{i-2}\dots a_{i-1}}),~~~i=2,\dots,l
\end{cases}$$
and $h_1=h_{\alpha_1}(t_1)\dots h_{\alpha_l}(t_l),~h_2=h_1^{-1}h$. By direct calculations it is easy to make sure that elements $h_1,h_2$ satisfy the necessary conditions.

\item $\Phi=B_l$. Cartan matrix of the root system $B_l$ can be transformed to the diagonal matrix $$diag(\underbrace{1,\dots,1}_{l-1}, 2)$$
      by integral elementary transformations of columns and rows. An arbitrary element of group  $H$ has the form
$$(t_1^2t_2^{-1},t_1^{-1}t_2^2t_3^{-1}, \dots,t_{l-3}^{-1}t_{l-2}^2t_{l-1}^{-1}, t_{l-2}^{-1}t_{l-1}^2t_l^{-2},t_{l-1}^{-1}t_l^2).$$
Let
$$\begin{cases}t_l=1,\\
t_{l-i}=({a_l^{i}a_{l-1}^{i-1}\dots a_{l-i+1}})^{-1},~~~i=1,\dots,l-1,
\end{cases}$$
then $h_1=h_{\alpha_1}(t_1)\dots h_{\alpha_l}(t_l),~h_2=h_1^{-1}h$ are required elements.
\item $\Phi=C_l$. By  integral elementary transformations of rows and columns Cartan matrix of the root system $C_l$ can be transformed to the diagonal matrix $$diag(\underbrace{1,\dots,1}_{l-1}, 2).$$
    An arbitrary element of $H$ has the form
$$(t_1^2t_2^{-1},t_1^{-1}t_2^2t_3^{-1},\dots, t_{l-2}^{-1}t_{l-1}^2t_{l}^{-1},t_{l-1}^{-2}t_l^2).$$
If
$$\begin{cases}t_1=1,\\
t_{i}=({a_1^{i-1}a_{2}^{i-2}\dots a_{i-1}})^{-1},~~~i=2,\dots,l,
\end{cases},$$
 then $h_1=h_{\alpha_1}(t_1)\dots h_{\alpha_l}(t_l),~h_2=h_1^{-1}h$ are required elements.

\item $\Phi=D_l$. By  integral elementary transformations Cartan matrix of the root system $D_l$ can be transformed to the following matrix
    \begin{align}
diag(\underbrace{1,\dots,1}_{l-1}, 4)~~~~~~~~~& l~\text{is odd}, \notag \\
diag(\underbrace{1,\dots,1}_{l-2}, 2,2)~~~~~~~& l~\text{is even}. \notag
\end{align}
The equation $T^4=a$ is solvable in the field $F$ for any $a$ if and only if the equations $T^2=a$ is solvable in $F$ for any  $a$.

An arbitrary element of the group $H$ is of the form
$$(t_1^2t_2^{-1},t_1^{-1}t_2^2t_3^{-1},\dots, t_{l-4}^{-1}t_{l-3}^2t_{l-2}^{-1},t_{l-3}^{-1}t_{l-2}^{2}t_{l-1}^{-1}t_l^{-1},t_{l-2}^{-1}t_{l-1}^2,t_{l-2}^{-1}t_l^2).$$
If
$$\begin{cases}t_1=t_l=1,\\
t_{i}=({a_1^{i-1}a_{2}^{i-2}\dots a_{i-1}})^{-1},~~~i=2,\dots,l-1,
\end{cases}$$
 then $h_1=h_{\alpha_1}(t_1)\dots h_{\alpha_l}(t_l),~h_2=h_1^{-1}h$ are required elements.

\item $\Phi=E_l~(l=6,7)$. By  integral elementary transformations of rows and columns Cartan matrix of the root system $E_l(l=6,7)$ can be transformed to the diagonal matrix
    \begin{align}
diag(\underbrace{1,\dots,1}_{l-1}, 3)~~~~~~~& l=6 \notag \\
diag(\underbrace{1,\dots,1}_{l-1}, 2)~~~~~~~& l=7 \notag
\end{align}

An arbitrary element of group $H$ has the following form
$$(t_1^2t_2^{-1},t_1^{-1}t_2^2t_3^{-1},\dots,t_{l-5}^{-1}t_{l-4}^2t_{l-3}^{-1},t_{l-4}^{-1}t_{l-3}^2t_{l-2}^{-1}t_{l-1}^{-1},t_{l-3}^{-1}t_{l-2}^2,t_{l-3}^{-1}t_{l-1}^2t_l^{-1},t_{l-1}^{-1}t_l^2).$$
Let
$$\begin{cases}t_1=t_{l-1}=1,\\
t_{i}=(a_1^{i-1}a_2^{i-2}\dots a_{i-1})^{-1},~~~i=2,\dots,l-2,\\
t_l=a_1^{l-4}a_2^{l-5}\dots a_{l-4}a_{l-1}^{-1},
\end{cases}$$
 and $h_1=h_{\alpha_1}(t_1)\dots h_{\alpha_l}(t_l),~h_2=h_1^{-1}h$, then $h_1, h_2$ are required elements.

\item $\Phi=E_8,F_4,G_2$. In this cases Cartan matrix can be transformed to identity matrix by integral elementary transformations of rows and columns. Therefore $H=\widehat{H}$.\hfill $\square$
\end{enumerate}
\begin{cor}\label{cor1} Let $G=\Phi(F)$ be an elementary Chevalley group of type $\Phi$ over the field $F$ and $h$ be an arbitrary element of group $\widehat{H}$. Then there exist such an elements $h_1 \in H,~h_2 \in \widehat{H}$, that $h=h_1h_2$ and in the Chevalley basis an element $h_2$ has the form
$$h_2=diag(\underbrace{1,\dots,1}_{|\Phi|-k},\underbrace{*,\dots,*}_k,\underbrace{1,\dots,1}_{|\Delta|}),$$
where $k$ is non-negative integer, which is given in the following table.
\begin{center}
\begin{tabular}[t]{|p{1em}|p{1em}|p{4em}|p{3em}|p{5em}|p{1em}|p{1em}|p{1em}|p{1em}|p{1em}|}
\hline
$\Phi$ & $A_l$ & ~~~~$B_l$ &~~~$C_l$&~~~~~$D_l$&$E_6$&$E_7$&$E_8$&$F_4$&$G_2$\\
\hline
k &$2l$&$2(2l-1)$&$l(l+1)$&$(l-1)(l+2)$&52&96&0&0&0 \\
\hline
\end{tabular}
\end{center}
\end{cor}
\textbf{Proof.}
In the Chevalley basis every element of group $\widehat{H}$ is obvious to have the diagonal form
 $$diag(\underbrace{*,\dots,*}_{|\Phi|},\underbrace{1,\dots,1}_{|\Delta|}).$$
 By the lemma \ref{lem20} there exist such an elements $h_1\in H$, $h_2\in \widehat{H}$ that $h=h_1h_2$ and $$h_2(x_{\alpha_i})=x_{\alpha_i},~~~i\in I.$$ The set $I$ doesn't depend on $h$ and it depends only on the type of root system $\Phi$. Let us denote $\overline{I}=\{1,2,\dots,l\}\setminus I$.

 It is obvious that if  $\alpha,~\beta,~\alpha+\beta$ are roots, such that $$h_2(x_{\alpha})=x_{\alpha}, ~~~h_2(x_{\beta})=x_{\beta},$$ then $h_2(x_{\alpha+\beta})=x_{\alpha+\beta}$. Therefore for any root $\alpha$ which is a linear combinations of simple roots from $\{\alpha_i~|~i\in I\}$ we have $h_2(x_{\alpha})=x_{\alpha}$. Hence $k$ is less or equal to the number of roots, which are such a linear combinations of simple roots, that at least one root from $\{\alpha_i~|~i\in \overline{I} \}$ has non-zero coefficient
By direct calculations we conclude that the number of such roots is equal to the number $k$ from the  condition of this corollary. ~\hfill $\square$

\textbf{{\scshape Example 1.}} Let us illustrate the last paragraph of the proof of corollary \ref{cor1} in the case of root system of type $A_l$. In this case we have  $I=\{1,\dots,l-1\},~\overline{I}=\{l\}$.

The set of roots, which are linear combinations of simple roots from  $\{x_i~|~i\in I\}=\{x_1,\dots,x_l\}$ forms a subsystem $A_{l-1}$ in the system $A_l$. Therefore the set of roots, which are such a linear combinations of simple roots, that the root $\alpha_l$ has non-zero coefficient, has the carinality $|A_l|-|A_{l-1}|=2l$.

\subsection{Rings and fields}

~~~~Let $F|L$ be an extension of the field $L$. Elements $x_1,\dots,x_k$ of the field $F$ are called \emph{algebraically independent} over $L$ if there is no such a polynomial $f(T_1,\dots,T_k)\neq0$, that the equality $f(x_1,\dots,x_k)=0$ holds.
Maximal set of elements of the field $F$, which are algebraically independent over the field $L$, is called transcendence basis of $F$ over $L$. Cardinality of transcendence basis of $F$ over $L$ doesn't depend on this basis and is called \emph{transcendence degree} of the field $F$ over $L$. The transcendence degree of $F$ over $L$ is denoted by ${\rm tr.deg}_LF$.

Let $\mathbb{Q}$ be the field of rational numbers, $\pi$ be the set of prime numbers and $2^{\pi}$ be the set of all subsets of $\pi$. Let us define the function $$\nu: \mathbb{Q}\rightarrow 2^{\pi}$$ by the rule: if $x=a/b\in \mathbb{Q}$, where $a$ and $b$ are mutually prime numbers, then $$\nu(x)=\{all~ prime~ devisos~of~ a\}\cup \{all~ prime~ devisors~of~ b\}.$$

\begin{lem}\label{lem1} Let $F$ be the field of zero characteristic, $x_1,\dots,x_k$ be elements of $F$, which are algebraically independent over $\mathbb{Q}$. Let $x_{k+1}$ be such an element of the field $F$, that elements $x_1,\dots, x_{k+1}$ are algebraically dependent over $\mathbb{Q}$. Let the automorphism $\delta$ of the field $F$ act on the elements $x_1,\dots, x_{k+1}$ by the following rule
$$
\delta: x_i\mapsto t_0 t_i x_i,~~~ i=1,\dots,k+1,
$$
where $t_0,\dots,t_{k+1}$ are such rational numbers, that $\nu(t_i)\cap\nu(t_j)=\varnothing$ for $i\neq j$ and $t_i\neq1$ for $i=1,\dots,k+1$. Then  $x_{k+1}=0$.
\end{lem}
\textbf{Proof.} Since elements $x_1,\dots, x_{k+1}$ are algebraically dependent over $\mathbb{Q}$, then we can choose such a polynomial of minimal degree $f(T_1,\dots,T_{k+1})\neq0$ with rational coefficients, that  $f(x_1,\dots, x_{k+1})=0$. Let
$$f(T_1,\dots,T_{k+1})=T_1^{n_1}\dots T_{k+1}^{n_{k+1}}+\sum_{i_1,\dots,i_{k+1}}a_{i_1,\dots,i_{k+1}}T_1^{i_1}\dots T_{k+1}^{i_{k+1}},$$
where  $T_1^{n_1}\dots T_{k+1}^{n_{k+1}}$ is a senior monomial.

If we define the map $\widetilde{\delta}:\{x_1,\dots,x_{k+1}\}\rightarrow\mathbb{Q}$ by the rule $\widetilde{\delta}(x_i)=t_0t_i$, then the equality $\delta(x_i)=\widetilde{\delta}(x_i)x_i$ is obvious.

Consider now the following polynomial $g(T_1,\dots,T_{k+1})$ with rational coefficients
$$g(T_1,\dots,T_{k+1})=f(\widetilde{\delta}(x_1)T_1,\dots,\widetilde{\delta}(x_{k+1})T_{k+1})$$
and note that $$g(x_1,\dots,x_{k+1})=f(\delta(x_1), \dots, \delta(x_{k+1}))=\delta(f(x_1,\dots, x_{k+1}))=\delta(0)=0.$$
Let
$$h(T_1,\dots,T_{k+1})=f(T_1,\dots,T_{k+1})-\widetilde{\delta}(x_1)^{-n_1} \dots \widetilde{\delta}(x_{k+1})^{-n_{k+1}}g(T_1,\dots,T_{k+1}).$$
Since polynomials $f$ and $g$ are equal to zero on the arguments $x_1,\dots,x_{k+1}$, then polynomial $h$  is also equal to zero on this arguments. Degree of the polynomial $h$ is less then the degree of $f$, and by the minimality of $f$ we conclude, that $h(T_1,\dots,T_{k+1})\equiv0$. Polynomial $h$ has the following form
$$h(T_1,\dots,T_{k+1})=\sum_{i_1,\dots,i_{k+1}}\left(1-\widetilde{\delta}(x_1)^{i_1-n_1} \dots \widetilde{\delta}(x_{k+1})^{i_{k+1}-n_{k+1}}\right)a_{i_1,\dots,i_{k+1}}T_1^{i_1}\dots T_{k+1}^{i_{k+1}}.$$
Since $h(T_1,\dots,T_{k+1})=0$, then
$$\left(1-\widetilde{\delta}(x_1)^{i_1-n_1} \dots \widetilde{\delta}(x_{k+1})^{i_{k+1}-n_{k+1}}\right)a_{i_1,\dots,i_{k+1}}=0$$
for all $i_1,\dots,i_{k+1}$. If $a_{i_1,\dots,i_{k+1}}\neq0$, then $1-\widetilde{\delta}(x_1)^{i_1-n_1} \dots \widetilde{\delta}(x_{k+1})^{i_{k+1}-n_{k+1}}=0$. Since $\widetilde{\delta}(x_i)=t_0t_i$, then we have
$$t_0^{(n_1+\dots +n_{k+1})-(i_1+\dots+i_{k+1})}
\prod _{j=1}^{k+1}t_j^{n_j-i_j}=1.$$

And therefore
$$\nu\left(t_0^{(n_1+\dots +n_{k+1})-(i_1+\dots+i_{k+1})}
\prod _{j=1}^{k+1}t_j^{n_j-i_j}\right)=\nu(1)=\varnothing,$$
but it contradicts to the conditions on $t_i$. Hence $a_{i_1,\dots,i_{k+1}}=0$ for all $i_1,\dots,i_{k+1}$. Therefore $f(T_1,\dots,T_{k+1})=T_1^{n_1}\dots T_{k+1}^{n_{k+1}}$ and from the equality $$0=f(x_1,\dots,x_{k+1})=x_1^{n_1}\dots x_{k+1}^{n_{k+1}}$$ follows, that $x_i=0$ for some $i$. Since $x_1,\dots,x_k$ are algebraically independent over $\mathbb{Q}$ then $x_{k+1}=0$. Lemma is proved.\hfill $\square$

Note that the condition $\nu(t_i)\cap\nu(t_j)=\varnothing$ for $i\neq j$ in the lemma \ref{lem1} can be replaced by weaker condition $t_0^{i_0}t_1^{i_1}\dots t_{k+1}^{i_{k+1}}\neq 1$ if $i_j\neq 0$ for some $j$.

The following statement is generalization of the lemma 1 from the paper \cite{Nas}.
\begin{lem}\label{lem3} Let $R$ be an integral domain and $M$ be an infinite subset of $R$. Let $f(T)$ be non-constant rational function with coefficients from the ring $R$. Then the set $P=\{f(a): a \in M\}$ is infinite.
\end{lem}
\textbf{Proof.} Let $f(T)=g(T)/h(T)$, where $g,h$ are polynomials of one variable $T$ with coefficients from the ring $R$. If we suppose that $P$ is finite, i. e. $P=\{b_1 \dots b_m\}$, then the polynomial $$p(T)=\prod_{i=1}^k(g(T)-b_ih(T))$$
has an infinite set of roots $M$, and therefore $p(T)\equiv0$. Hence $g(T)-b_ih(T)\equiv0$ for some $i$, i. e. $f(T)\equiv b_i$, but it contradicts to the conditions of lemma.\hfill $\square$

\subsection{Automorphisms of Chevalley groups}
Let us remind some classical types of automorphisms of elementary Chevalley group $G=\Phi(F)$ (\cite[\S 12.3]{Car}).

\emph{Inner automorphisms.} For an arbitrary element $h\in G$ symbol $\varphi_h$ means an automorphism, which acts by the following way
$$\varphi_h: g\rightarrow hgh^{-1}.$$
Automorphism $\varphi_h$ is called an inner automorphism of group $G$ induced by the element $h$. all the inner automorphisms form a normal subgroup in the group of all automorphisms of $G$.

 \emph{Diagonal automorphisms.}
As we already noted, the group $\widehat{H}$ normalizes elementary Chevalley group $G$ in the group of all automorphisms of Lie algebra $\mathcal{L}$. Hence if $h \in \widehat{H}$ then the map
$$\varphi_h: g \mapsto hgh^{-1}$$
is an automorphism of group $G$. If $h$ belongs to $\widehat{H}$ and doesn't belong to $H$ then this automorphism is called \emph{diagonal automorphism}. If $h\in H$, then $\varphi_h$ is an inner automorphism. In the Chevalley basis all the diagonal automorphisms are induced by conjugations by diagonal matrices.

\emph{Field automorphisms.}
If $\delta$ is an automorphism of the field $F$, then the map
$$\overline{\delta}: x_{\alpha}(t) \mapsto x_{\alpha}(\delta(t)),~~~\alpha \in \Phi, t\in F $$
can be extended to the automorphism of group $G$. Such an automorphism is called \emph{field automorphism}. Automorphism $\overline{\delta}$ maps the element $h_{\alpha}(t)$ to the element $h_{\alpha}(\delta(t))$.

\emph{Graph automorphisms.}
Automorphisms of this type arise from symmetries of Dynkin diagram. Symmetry of Dynkin diagram of the root system $\Phi$ is a permutation $\rho$ of nodes of Dynkin diagram, such that the number of edges between the nodes $\alpha_i$ and $\alpha_j$ is equal to the number of edges between $\rho(\alpha_i)$ and $\rho(\alpha_j)$ for any pair $i\neq j$.

If all the simple roots has the same length, then any permutation of simple roots induces permutation of all roots, which is denote by the same symbol $\rho$. Then there exist such numbers $\gamma_{\alpha}=\pm1$, that the map

$$\overline{\rho}:  x_{\alpha}(t) \mapsto x_{\rho(\alpha)}(\gamma_{\alpha}t),~~~\alpha \in \Phi, t\in F$$
can be extended to the automorphism of $G$. It is possible to choose such a numbers $\gamma_{\alpha}$, that $\gamma_{\alpha}=1$ for $\alpha \in \Delta$ or $-\alpha \in \Delta$. Graph automorphism maps $h_{\alpha}(t)$ to $h_{\rho(\alpha)}(t)$.
If characteristic of the field $F$ is equal to zero, then graph automorphisms exist only in the cases when all the simple roots have the same length (i. e. for root systems of types $A_l, D_l, E_6$). Symmetries of Dynkin diagrams are the following
\begin{center}
\includegraphics[width=65mm]{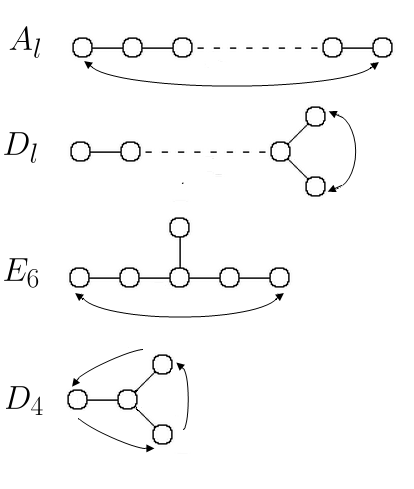}
\end{center}
Order of graph automorphism is equal to $2$ or $3$.

 In \cite{STR, JHam} it is proved, that for any automorphism  $\varphi$ of elementary Chevalley group $G=\Phi(F)$ of type $\Phi$ over the field $F$ there exist inner automorphism $\varphi_g$, diagonal automorphism $\varphi_h$, graph automorphism $\overline{\rho}$ and field automorphism $\overline{\delta}$, such that $\varphi=\overline{\rho}\overline{\delta}\varphi_h\varphi_g$.

Note that graph automorphisms and field automorphisms commute.
\section{Proofs of the main results}
\begin{ttt}\label{ttt1} Let  $G$ be a Chevalley group of type $\Phi\neq A_1$ over the field $F$ of zero characteristic and transcendence degree of $F$ over $\mathbb{Q}$ is finite. Then
\begin{enumerate}
\item If $\Phi$ has one of types $A_l(l\geq 7)$, $B_l(l\geq 4)$, $E_8,~F_4,~G_2$, then $G$ possess $R_{\infty}$ property.
\item If the equation $f(T)=a$ is solvable in $F$ for any $a$, where $f$ is a polynomial from the lemma \ref{lem20}, then $G$ also posses $R_{\infty}$ property in the cases of root systems of  types $A_l(l=2,3,4,5,6),~B_l(l=2,3),~C_l(l\geq 3),~D_l(l\geq 4),~E_6,~E_7$.
\end{enumerate}
\end{ttt}
\textbf{Proof.} Since $G/Z(G)\cong \Phi(F)$ then, according to proposition \ref{pr3} it is sufficient to prove theorem for  $G=\Phi(F)$.

 Consider the automorphism $\varphi$ of group $G$ and show that $R(\varphi)=\infty$.
Let $\varphi\in {\rm Aut}~G$ then $\varphi=\overline{\rho}\overline{\delta}\varphi_h\varphi_g$ for suitable inner, diagonal, graph and field automorphisms $\varphi_g,~\varphi_h,~\overline{\rho},~\overline{\delta}$ respectively.  By the corollary \ref{cor1} for an element $h$ we have $h=h_1h_2$, where $h_1\in H$, $h_2\in \widehat{H}$ and in Chevalley basis an element $h_2$ has the form
$$h_2=diag(\underbrace{1,\dots,1}_{|\Phi|-k},\underbrace{*,\dots,*}_k,\underbrace{1,\dots,1}_{|\Delta|})$$
for suitable $k$. Therefore $\varphi_h=\varphi_{h_1h_2}=\varphi_{h_1}\varphi_{h_2}=\varphi_{h_2}\varphi_{h_1}$ and automorphism $\varphi$ has the following form
$$\varphi=\overline{\rho}\overline{\delta}\varphi_h\varphi_g=\overline{\rho}\overline{\delta}\varphi_{h_2}\varphi_{h_1}\varphi_g=\overline{\rho}\overline{\delta}\varphi_{h_2}\varphi_{h_1g}.$$
By the proposition \ref{pr1} Reidemeister number $R(\varphi)$ is infinite if and only if the number $R(\varphi\varphi_{g^{-1}h_1^{-1}})$ is infinite, hence it is sufficient to prove the statement for $\varphi=\overline{\rho}\overline{\delta}\varphi_{h_2}$.

Suppose that $R(\varphi)<\infty$ and consider the following elements of group $G$
 $$g_i=h_{\alpha_1}(p_{i1})h_{\alpha_2}(p_{i2})\dots{h_{\alpha_l}}(p_{il}),~~~ i=1,2,\dots,$$
 where $p_{11}<p_{12}<\dots<p_{1l}<p_{21}<p_{22}<\dots$ are prime integers. In the Chevalley basis an element $g_i$ has the form
 $$g_i=diag(a_{i1},a_{i2},\dots,a_{i|\Phi|},\underbrace{1,\dots,1}_{|\Delta|}),$$
for rational numbers $a_{ij}$, such that $\nu(a_{ij})\neq\varnothing$ and $\nu(a_{ij})\cap\nu(a_{rs})=\varnothing$ for $i\neq r$ (since $\nu(a_{ij})\subseteq\{p_{i1},p_{i2},\dots,p_{il}\}$).

Since the number of $\varphi$-conjugacy classes is finite, then there exist an infinite subset of  $\varphi$-conjugated elements in  $g_1,g_2,\dots$ . Without loosing of generality we can consider that all the elements $g_1,g_2,\dots$ belong to $[g_1]_{\varphi}$. Then there exist such an elements $Z_1,Z_2,\dots$ in $G$, that the following equalities hold
$$
g_i=Z_ig_1\varphi(Z_i^{-1}),~~~ i=2,3,\dots
$$
Since $\varphi=\overline{\rho}\overline{\delta}\varphi_{h_2}$ then we have \begin{equation}\label{eq0}
g_i=Z_ig_1\overline{\rho}\overline{\delta}\varphi_{h_2}(Z_i^{-1}),~~~ i=2,3,\dots
\end{equation}
and depending on the type of root system $\Phi$ we have 4 cases.

\emph{Case 1.} $\Phi$ is a root system of type $E_8,~F_4,~G_2$, or $\Phi$ is a root system of type $A_2,~A_3,~B_2,~B_3,~C_l(l\geq 3),~D_l(l\geq 5),~E_6,~E_7$ and equation $f(T)=a$ is solvable in $F$ for any $a$, where $f$ is a polynomial from the lemma \ref{lem20}. In this case we have $\overline{\rho}^2=id$, $\varphi_{h_2}=id$, and therefore $\varphi=\overline{\rho}\overline{\delta}$. Then the equality (\ref{eq0}) can be rewritten by the following way
 \begin{equation}\label{eq1}
 g_i=Z_ig_1\overline{\rho}\overline{\delta}(Z_i^{-1}),~~~i=2,3,\dots
 \end{equation}
Acting on this equality by the automorphism $\varphi=\overline{\rho}\overline{\delta}$ considering the fact that graph automorphism and field automorphism commute, we have
\begin{equation}\label{eq2}
\overline{\rho}\overline{\delta}(g_i)=\overline{\rho}\overline{\delta}(Z_i)\overline{\rho}\overline{\delta}(g_1)\overline{\rho}^2\overline{\delta}^2(Z_i^{-1})=\overline{\rho}\overline{\delta}(Z_i)\overline{\rho}\overline{\delta}(g_1)\overline{\delta}^2(Z_i^{-1}),~~~ i=2,3,\dots
\end{equation}
If we multiply the equalities (\ref{eq1}) and (\ref{eq2}), then we have
\begin{equation}\label{eq3}
g_i\overline{\rho}\overline{\delta}(g_i)=Z_ig_1\overline{\rho}\overline{\delta}(g_1)\overline{\delta}^2(Z_i^{-1}),~~~ i=2,3,\dots
\end{equation}
 Since $\delta$ acts identically on the prime subfield $\mathbb{Q}$ of the field $F$ and $\overline{\delta}(h_{\alpha_i}(t))=h_{\alpha_i}(\delta(t))$, then $\overline{\delta}(h_{\alpha_i}(t))=h_{\alpha_i}(t)$ for every $t\in\mathbb{Q}$ and therefore $\overline{\delta}(g_i)=g_i$ and
$$\overline{\rho}\overline{\delta}(g_i)=\overline{\rho}(g_i),~~~i=2,3,\dots$$
If we denote $A_i=g_i\overline{\rho}(g_i)$ then in Chevalley basis matrix $A_i$ has the following diagonal form
$$A_i=diag(b_{i1},b_{i2},\dots,b_{i|\Phi|},\underbrace{1,\dots,1}_{|\Delta|}),$$
where $\nu(b_{ij})\neq\varnothing$ and $\nu(b_{ij})\cap\nu(b_{rs})=\varnothing$ for $i\neq r$ (since $\nu(b_{ij})\subseteq\underset{j}{\bigcup} \nu(a_{ij})$).

If we denote $\overline{\theta}=\overline{\delta}^2$ then $\overline{\theta}$ is also a field automorphism. In this denotation the equality (\ref{eq3}) can be rewritten by
$A_i=Z_iA_1\overline{\theta}(Z_i^{-1})$ and therefore
\begin{equation}\label{eq4}
\overline{\theta}(Z_i)=A_i^{-1}Z_iA_1.
\end{equation}
Let the matrix $Z_i$ has the form $Z_i=(z_{i,rs})_{r,s=1}^{|\Phi|+|\Delta|}$ in the Chevalley basis. Then the equality (\ref{eq4}) can be rewritten in more details
$$
\begin{vmatrix}
\theta(z_{i,11})&\dots\\
\vdots&\ddots\\
\end{vmatrix}
=
\begin{vmatrix}
b_{i1}^{-1}{b_{11}}z_{i,11}&\dots\\
\vdots&\ddots\\
\end{vmatrix}.
$$
Then we have the following system of equalities
$$\begin{cases}
\theta(z_{i,rs})=b_{ir}^{-1}b_{1s}z_{i,rs},~~~ r,s=1,\dots,|\Phi|,\\
\theta(z_{i,rs})=b_{ir}^{-1}z_{i,rs},~~~ r=1,\dots,|\Phi|, s=|\Phi|+1,\dots,|\Phi|+|\Delta|
\end{cases},~~~i=2,3,\dots$$

In the set $z_{i,11}$ (i=1,2,\dots) we can choose a maximal set of algebraically independent elements. This set is finite since the transcendence  degree of $F$ over $\mathbb{Q}$ is finite. Without loosing of generality we can consider that the first $t$ elements $z_{1,11},z_{2,11},\dots,z_{t,11}$ are algebraically independent over $\mathbb{Q}$. Similarly find maximal sets of algebraically indepent elements $\left\{z_{i_{j},rs}\right\}_{j=1}^{n_{rs}}$ in the sets ${z_{i,rs}}$ ($r=1,\dots,|\Phi|,~s=1,\dots,|\Phi|+|\Delta|$) and denote by $n$ the following value $$n=\underset{r,s}{{\rm max}}(t,i_{n_{rs}}).$$
By the definition, the number $n$ is such a number, that elements $z_{1,rs},\dots,z_{n,rs},z_{(n+1),rs}$ are algebraically dependent over $\mathbb{Q}$ for any  $r=1,\dots,|\Phi|,~s=1,\dots,|\Phi|+|\Delta|$.

Let $t_0=b_{11}, t_i=b_{i1}^{-1}$, $i=2,3,\dots$, then the following statements hold
\begin{enumerate}
\item Elements $z_{1,11},z_{2,11},\dots,z_{t,11}$ are algebraically independent over $\mathbb{Q}$.
\item Elements $z_{1,11},z_{2,11},\dots,z_{t,11}, z_{(n+1),11}$ are algebraically dependent over $\mathbb{Q}$.
\item $\nu(t_i)\neq \varnothing$ for $i\neq0$.
\item $\nu(t_i) \cap \nu(t_j)=\varnothing$ for $i\neq j$.
\end{enumerate}
By the lemma \ref{lem1} we conclude that $z_{(n+1),11}=0$. By the similar arguments for every $z_{(n+1),rs}$ for $r=1,\dots,|\Phi|,~s=1,\dots,|\Phi|+|\Delta|$ we have that
$$Z_{n+1}=
\begin{vmatrix}
O_{|\Phi|\times(|\Phi|+l)}\\
*&
\end{vmatrix},
$$
but it contradicts to the fact, that $Z_{n+1}$ is invertible. Therefore all the matrices $A_1,A_2,\dots$ can't be $\varphi$-conjugated, and hence $R(\varphi)=\infty$.

\emph{Case 2.} The root system $\Phi$ is of the type $B_l~(l\geq 4)$. In this case $|\Delta|=l$, $|\Phi|=2l^2$, $k=2(2l-1)$ and $\overline{\rho}=id$. Then the equality (\ref{eq0}) has the form

\begin{equation}\label{eq5}
 g_i=Z_ig_1\overline{\delta}\varphi_{h_2}(Z_i^{-1}),~~~i=2,3,\dots
 \end{equation}
 Let the matrix $Z_i$ be splitted on the blocks
 $$\newcommand{\tempca}{\multicolumn{1}{c|}{Z_{i,11}}}
\newcommand{\tempcaa}{\multicolumn{1}{c|}{Z_{i,12}}}
\newcommand{\tempcab}{\multicolumn{1}{c|}{Z_{i,21}}}
\newcommand{\tempcac}{\multicolumn{1}{c|}{Z_{i,22}}}
\newcommand{\tempcad}{\multicolumn{1}{c|}{Z_{i,31}}}
\newcommand{\tempcae}{\multicolumn{1}{c|}{Z_{i,32}}}
\newcommand{\tempcag}{\multicolumn{1}{c|}{Z_{i,13}}}
\newcommand{\tempcagg}{\multicolumn{1}{c|}{Z_{i,23}}}
\newcommand{\tempcaggg}{\multicolumn{1}{c|}{Z_{i,33}}}
\newcommand{\tempcTT}{\multicolumn{1}{c|}{{|\Phi|-k}}}
\newcommand{\tempcTTT}{\multicolumn{1}{c|}{k}}
\newcommand{\tempcTTTT}{\multicolumn{1}{c|}{l}}
\begin{matrix}
&|\Phi|-k&k&l\\\cline{2-4}
\tempcTT& \tempca&\tempcaa&\tempcag\\\cline{2-4}
\tempcTTT&\tempcab&\tempcac&\tempcagg\\\cline{2-4}
\tempcTTTT&\tempcad&\tempcae&\tempcaggg\\\cline{2-4}
&&&
\end{matrix}~~,
$$
 where the number of rows of block is written on the left of matrix and the number of columns is written on the top of matrix. Then the matrix $\varphi_{h_2}(Z_i)$ has the form
$$\newcommand{\tempca}{\multicolumn{1}{c|}{Z_{i,11}}}
\newcommand{\tempcaa}{\multicolumn{1}{c|}{*}}
\newcommand{\tempcab}{\multicolumn{1}{c|}{*}}
\newcommand{\tempcac}{\multicolumn{1}{c|}{*}}
\newcommand{\tempcad}{\multicolumn{1}{c|}{Z_{i,31}}}
\newcommand{\tempcae}{\multicolumn{1}{c|}{*}}
\newcommand{\tempcag}{\multicolumn{1}{c|}{Z_{i,13}}}
\newcommand{\tempcagg}{\multicolumn{1}{c|}{*}}
\newcommand{\tempcaggg}{\multicolumn{1}{c|}{Z_{i,33}}}
\newcommand{\tempcTT}{\multicolumn{1}{c|}{{|\Phi|-k}}}
\newcommand{\tempcTTT}{\multicolumn{1}{c|}{k}}
\newcommand{\tempcTTTT}{\multicolumn{1}{c|}{l}}
\varphi_{h_2}(Z_i)=h_2Z_ih_2^{-1}=~~\begin{matrix}
&|\Phi|-k&k&l\\\cline{2-4}
\tempcTT& \tempca&\tempcaa&\tempcag\\\cline{2-4}
\tempcTTT&\tempcab&\tempcac&\tempcagg\\\cline{2-4}
\tempcTTTT&\tempcad&\tempcae&\tempcaggg\\\cline{2-4}
&&&
\end{matrix}~~.$$
From the equality (\ref{eq5}) we have the following equality

\begin{equation}
 \overline{\delta}\varphi_{h_2}(Z_i)=g_i^{-1}Z_ig_1,~~~i=1,2,\dots
 \end{equation}
Then by  similar to the case 1 arguments we conclude that there exist such a number $n$, that the matrix $Z_n$ is of the form
$$  \newcommand{\tempca}{\multicolumn{1}{c|}{O_{(|\Phi|-k)\times(|\Phi|-k)}}}
\newcommand{\tempcaa}{\multicolumn{1}{c|}{Z_{i,12}}}
\newcommand{\tempcab}{\multicolumn{1}{c|}{Z_{i,21}}}
\newcommand{\tempcac}{\multicolumn{1}{c|}{Z_{i,22}}}
\newcommand{\tempcad}{\multicolumn{1}{c|}{Z_{i,31}}}
\newcommand{\tempcae}{\multicolumn{1}{c|}{Z_{i,32}}}
\newcommand{\tempcag}{\multicolumn{1}{c|}{O_{(|\Phi|-k)\times l}}}
\newcommand{\tempcagg}{\multicolumn{1}{c|}{Z_{i,23}}}
\newcommand{\tempcaggg}{\multicolumn{1}{c|}{Z_{i,33}}}
\newcommand{\tempcTT}{\multicolumn{1}{c|}{{|\Phi|-k}}}
\newcommand{\tempcTTT}{\multicolumn{1}{c|}{k}}
\newcommand{\tempcTTTT}{\multicolumn{1}{c|}{l}}
\begin{matrix}
&|\Phi|-k&k&l\\\cline{2-4}
\tempcTT& \tempca&\tempcaa&\tempcag\\\cline{2-4}
\tempcTTT&\tempcab&\tempcac&\tempcagg\\\cline{2-4}
\tempcTTTT&\tempcad&\tempcae&\tempcaggg\\\cline{2-4}
&&&
\end{matrix}~~.
$$

 Since $|\Phi|=2l^2$ and $k=2(2l-1)$, then for $l\geq 4$ inequality $|\Phi|-k>k$ holds, and therefore the first $|\Phi|-k$ rows of matrix $Z_n$ are linear dependent, but it contradicts to the non-singularity of matrix $Z_n$.

\emph{Case 3.} The root system $\Phi$ has the type $A_l~(l\geq 3)$. In this case $|\Delta|=l$, $|\Phi|=l(l+1)$, $ k=2l$, $\overline{\rho}^2=id$ and $\varphi=\overline{\rho}\overline{\delta}\varphi_{h_2}$. If we act by the automorphism $\varphi$ on the equality (\ref{eq0}), then we have

\begin{equation}\label{eq8}
 \varphi(g_i)=\varphi(Z_i)\varphi(g_1)\varphi^2(Z_i^{-1}),~~~i=1,2,\dots
 \end{equation}
After multiplication of the equalities (\ref{eq0}) and (\ref{eq8}) we have the following equality
\begin{equation}\label{eq9}
 g_i\varphi(g_i)=Z_ig_1\varphi(g_1)\varphi^2(Z_i^{-1}),~~~i=1,2,\dots
 \end{equation}
In the case 1, we have already noted, that  $\overline{\delta}(g_i)=g_i$. Moreover, since $h_2$ and $g_i$ has the diagonal form in the Chevalley basis, then $\varphi_{h_2}(g_i)=g_i$. Therefore $\varphi(g_i)=\overline{\rho}(g_i)$. Now we denote $A_i=g_i\overline{\rho}(g_i)$ and rewrite the equality (\ref{eq8}) using the fact, that $\varphi=\overline{\rho}\overline{\delta}\varphi_{h_2}:$

 \begin{equation}\label{eq10}
 A_i=Z_iA_1\overline{\rho}\overline{\delta}\varphi_{h_2}\overline{\rho}\overline{\delta}\varphi_{h_2}(Z_i^{-1}),~~~i=2,3,\dots
 \end{equation}

For an arbitrary  $h(\chi)\in \widehat{H}$ the following identities hold
 $$\varphi_{h(\chi)}\overline{\delta}=\overline{\delta}\varphi_{h(\delta^{-1}\circ\chi)},$$
$$\varphi_{h(\chi)}\overline{\rho}=\overline{\rho}\varphi_{h(\chi_1)},$$
where $\chi_1(\alpha)=\chi(\rho^{-1}(\alpha))$ for $\alpha\in \Phi$. Since $h_2\in \widehat{H}$, then $h_2=h(\chi)$ for suitable $\chi$ and hence $$\overline{\rho}\overline{\delta}\varphi_{h(\chi)}\overline{\rho}\overline{\delta}\varphi_{h(\chi)}=
\overline{\rho}\overline{\delta}\overline{\rho}\varphi_{h(\chi_1)}\overline{\delta}\varphi_{h(\chi)}=
\overline{\rho}\overline{\delta}\overline{\rho}\overline{\delta}\varphi_{h(\delta^{-1}\circ\chi_1)}\varphi_{h(\chi)}=
\overline{\delta}^2\varphi_{h(\delta^{-1}\circ\chi_1)h(\chi)}.$$

For the root system $\Phi$ of type $A_l$ permutation $\rho$ acts on the simple root $\alpha_i$ by the way $\rho(\alpha_i)=\alpha_{l-i+1}$. Hence $\rho(\alpha_1+\dots+\alpha_l)=\alpha_1+\dots+\alpha_l$, and all other roots, which are linear combinations of simple roots witf non-zero coefficient of $\alpha_1$, are mapped under the map $\rho$ to the roots, which have zero coefficient of $\alpha_1$ in its presentation. Therefore in the Chevalley basis the element $h(\delta^{-1}\circ\chi_1)$ has the form
$$h(\delta^{-1}\circ\chi_1)=diag(\underbrace{1,\dots,1}_{|\Phi|-2k+1},\underbrace{*,\dots,*}_{k-1},\underbrace{1,\dots,1}_{k-1},*,\underbrace{1,\dots,1}_{|\Delta|}).$$
then the element $h_3=h(\delta^{-1}\circ\chi_1)h(\chi)$ has the following form in Chevalley basis
$$h_3=diag(\underbrace{1,\dots,1}_{|\Phi|-2k+1},\underbrace{*,\dots,*}_{2k-1},\underbrace{1,\dots,1}_{|\Delta|}).$$
Therefore the equality (\ref{eq10}) can be rewritten
$$
 A_i=Z_iA_1\overline{\delta}^2\varphi_{h_3}(Z_i^{-1}),~~~i=1,2,\dots
$$
By the similar to the case 2 arguments, using the fact that for the root system of type $A_l$ ($l\geq 7$) the inequality $|\Phi|-2k+1>2k-1$ holds, we conclude that the matrix $Z_n$ is singular for big enough number $n$. It contradicts the fact that $Z_n$ belongs to $G$ and therefore must be invertible.

\emph{Case 4.} Root system $\Phi$ has type $D_4$. In this case we have $\overline{\rho}^2=id$, either $\overline{\rho}^3=id$. If equation $T^2=a$ is solvable in the field $F$ for any $a$, then $\widehat{H}=H$, therefore $\varphi_{h_2}=id$ and $\varphi=\overline{\rho}\overline{\delta}$.

The case $\overline{\rho}^2=id$ is similar to the case 1. If $\overline{\rho}^3=id$, then acting on the equality (\ref{eq0}) by the automorphism $\varphi$ we have
\begin{eqnarray}
\nonumber  g_i&=&Z_ig_1\varphi(Z_i^{-1}) ,\\
\nonumber  \varphi(g_i)&=&\varphi(Z_i)\varphi(g_1)\varphi^2(Z_i^{-1}),\\
\nonumber \varphi^2(g_i)&=&\varphi^2(Z_i)\varphi^2(g_1)\varphi^3(Z_i^{-1}).
\end{eqnarray}
If we multiply this equalities, we have
$$g_i\varphi(g_i)\varphi^2(g_i)=Z_ig_1\varphi(g_1)\varphi^2(g_1)\varphi^3(Z_i^{-1}).$$
After denoting $A_i=g_i\varphi(g_i)\varphi^2(g_i)$ using the fact that  $\varphi^3=\overline{\rho}^3\overline{\delta}^3=\overline{\delta}^3$, we can apply arguments of case 1.

Since $\varphi$ is an arbitrary automorphism of group $G$, we conclude that $G$ possess $R_{\infty}$ property.\hfill $\square$

Consider the field $F(T)$ of rational functions of one variable $T$ with coefficients from the field $F$ and let  $\Phi(F(T))$ be elementary Chevalley group over this field. Elements of this group are matrix with the entries from $F(T)$. We can speak about value of rational function on the arbitrary element of the field  $F$ if its denominator is not equal to zero.

 Let $a\in F$. If all the entries of the elements $h(T)\in\Phi(F(T))$ are determined on the element $a$, then $h(a)\in \Phi(F)$.

The following technical lemma is extremely useful in the proof of the theorem 2.
\begin{lem} \label{lem7}Let $g(T)=h_{\alpha_1}(T)h_{\alpha_2}(T)\dots h_{\alpha_l}(T)$ be an element of group $\Phi(F(T))$ and $\chi:\mathbb{Z}\Phi\rightarrow F^*$ be a homomorphism from the additive group $\mathbb{Z}\Phi$ to the multiplicative group of the field $F$. Then for any
 $m$  an element  $g(T)^mh(\chi)$ in the Chevalley basis has such a diagonal form, that its trace belongs to $F(T)\setminus F$.
\end{lem}
\textbf{Proof.} It is obvious that in Chevalley basis an element $g(T)$ has the diagonal form with degrees of variable $T$ on the diagonal. Then trace of the element $g(T)^mh(\chi)$ is a rational function.

For showing that this rational function is not a constant it is sufficient to show that the absolute value of degree of one of diagonal elements is greater than  absolute value of degrees of all other diagonal elements.

Below we have written such a root, which corresponds to the maximal absolute value of variable $T$ depending on the root system $\Phi$:
\begin{enumerate}
\item $\Phi=A_l$,~~~$\alpha=\alpha_1+\alpha_2+\dots+\alpha_l$;
\item $\Phi=B_l$,~~~$\alpha=\alpha_1+\alpha_2+\dots+\alpha_{l-1}+2\alpha_l$;
\item $\Phi=C_l$,~~~$\alpha=2\alpha_1+2\alpha_2+\dots+2\alpha_{l-1}+\alpha_l$;
\item $\Phi=D_l$,~~~$\alpha=\alpha_1+\alpha_2+\dots+\alpha_l$;
\item $\Phi=E_l$,~~~$\alpha=\alpha_1+\alpha_2+\dots+\alpha_l$;
\item $\Phi=F_4$,~~~$\alpha=-2\alpha_2-\alpha_3$;
\item $\Phi=G_2$,~~~$\alpha=-3\alpha_1-\alpha_2$.
\end{enumerate}
Therefore for any root system we have constructed required element, and hence trace of the element $g(T)^mh(\chi)$ is non-constant rational function.\hfill $\square$

\textbf{{\scshape Example 2.}} Let us illustrate the proof of the lemma \ref{lem7} on the root system of type $A_l$. In this case the system of positive roots has the form
$$\Phi^+=\{\alpha_i+\alpha_{i+1}+\dots+\alpha_j~|~1\leq i < j\leq l\}\cup \Delta.$$
By the proposition \ref{pr2}, an element $g(T)$ acts on the elements $x_{\alpha}$, which correspond to simple roots, by the following way
$$g(T)(x_{\alpha_i})=\begin{cases} Tx_{\alpha_i},~~~i=1,l,\\
x_{\alpha_i},~~~i=2,\dots,l-1.
\end{cases}$$
Therefore we have
$$g(T)(x_{\alpha_i+\alpha_{i+1}+\dots+\alpha_j})=\begin{cases} Tx_{\alpha_i+\alpha_{i+1}+\dots+\alpha_j},~~~1=i<j<l,\\
Tx_{\alpha_i+\alpha_{i+1}+\dots+\alpha_j},~~~1<i<j=l,\\
x_{\alpha_i+\alpha_{i+1}+\dots+\alpha_j},~~~1<i<j<l,\\
T^2x_{\alpha_i+\alpha_{i+1}+\dots+\alpha_j},~~~i=1,~j=l.
\end{cases}$$
Hence the maximal degree of variable $T$, which can be found on the diagonal of matrix $g(T)$, corresponds to the root $\alpha_1+\alpha_{i+1}+\dots+\alpha_l$ and is equal to 2.

Denote the trace of matrix $g(T)^mh$ by $\psi_{h,m}(T)$.

\begin{ttt}\label{ttt2} Let  $G$ be a Chevalley group of type $\Phi\neq A_1$ over the field $F$ of zero characteristic. If an automorphism group of the field $F$ is torsion, then $G$ possess $R_{\infty}$ property.
\end{ttt}
\textbf{Proof.} Since $G/Z(G)\cong \Phi(F)$ then by proposition \ref{pr3} we can consider  $G=\Phi(F)$.
Let $\varphi\in {\rm Aut}~G$, then $\varphi=\overline{\rho}\overline{\delta}\varphi_h\varphi_g$ for suitable inner, diagonal, graph and field automorphisms $\varphi_g,~\varphi_h,~\overline{\rho},~\overline{\delta}$ respectively.
Because of proposition \ref{pr1} we can consider that $\varphi=\overline{\rho}\overline{\delta}\varphi_{h}$. Since an automorphism group of the field $F$ is torsion, then any automorphism $\overline{\delta}$ has finite order $m$.

Let $g(T)$ be an element of group $\Phi(F(T))$, which is defined in lemma \ref{lem7}. Suppose that $R(\varphi)<\infty$ and consider the set of elements $g_i=g(x_i)$, where $\{x_i\}$ is an infinite set of non-zero rational numbers. Without loosing of generality we can consider that all the elements $g_i$ belong to the class $[g_1]_{\varphi}$. Hence there exist such an elements $Z_1,Z_2,\dots,$ that the following equalities hold
$$
g_i=Z_ig_1\varphi(Z_i^{-1}),~~~~i=2,3,\dots
$$
 Acting on this equalities bydegrees of the automorphism $\varphi$ we have
\begin{eqnarray}\label{eq12}
\nonumber  g_i&=&Z_ig_1\varphi(Z_i^{-1}) ,\\
\nonumber  \varphi(g_i)&=&\varphi(Z_i)\varphi(g_1)\varphi^2(Z_i^{-1}),\\
 \varphi^2(g_i)&=&\varphi^2(Z_i)\varphi^2(g_1)\varphi^3(Z_i^{-1}),\\
\nonumber &\dots&\\
\nonumber \varphi^{m-1}(g_i)&=&\varphi^{m-1}(Z_i)\varphi^{m-1}(g_1)\varphi^{m}(Z_i^{-1}).
\end{eqnarray}
Since graph automorphisms commute with field automorphism, and all the diagonal automorphisms form a normal subgroup in the group, which is generated by all the diagonal, graph and field automorphisms, then we have
$$\varphi^i=(\overline{\rho}\overline{\delta}\varphi_{h})^i=\overline{\rho}^i\overline{\delta}^i\varphi_{h^{\prime}}.$$
Since in the Chevalley basis an element $g_i$ has the diagonal form, and diagonal automorphism acts by conjugation by diagonal matrix, then $\varphi_h(g_i)=g_i$ for every diagonal automorphism $\varphi_h$. Moreover since any automorphism of the field $F$ acts identically on prime subfield, then $\overline{\delta}({g_i})=g_i$. Also for any graph automorphism  $\overline{\rho}$ we have $\overline{\rho}(g_i)=\overline{\rho}(h_{\alpha_1}(x_i)\dots h_{\alpha_l}(x_i))=h_{\rho(\alpha_1)}(x_i)\dots h_{\rho(\alpha_l)}(x_i)=g_i$. Hence $\varphi(g_i)=g_i$ for any $i$. Using this notes multiply all the equalities from (\ref{eq12})
$$g_i^m=Z_ig_1^m\varphi_{h^{\prime}}(Z_i^{-1}).$$
After multiplying this equality on $h^{\prime}$ from the right we have
$$g_i^mh^{\prime}=Z_ig_1^mh^{\prime}Z_i^{-1},$$
i.e. $g_i^mh^{\prime}$ and $g_1^mh^{\prime}$ are conjugated, and therefore their traces coincide, i.e. $\psi_{h^{\prime},m}(x_i)=\psi_{h^{\prime},m}(x_1)$ for any $i=2,3,\dots,$ but it contradicts to lemma \ref{lem3}.\hfill $\square$

\section{Twisted conjugacy class of the unit element}
In this section we state some sufficient conditions   when the twisted   conjugacy class of the unit element is a subgroup of Chevalley group. At first, let us remind three propositions from the paper \cite{BarNasNes}.

\begin{prp}\label{pr4} Let $G$ be a group and $\varphi$ be a central automorphism of group $G$. Then the set $[e]_{\varphi}$ is a subgroup of $G$.
\end{prp}

\begin{prp}\label{pr5} If  $\varphi$-conjugacy class $[e]_{\varphi}$ of unit element $e$ of group $G$ is a subgroup of $G$, then this subgroup is normal subgroup of $G$.
\end{prp}

If $N$ is a normal $\varphi$-admissible subgroup of $G$, then denote by $\overline{e}$ a unit element of quotient group $G/N$, and by $\overline{\varphi}$ an automorphism, which is induced by $\varphi$ on the quotient group.
\begin{prp}\label{pr6} Let $G$ be such a group, that for any automorphism $\varphi$ of this group $G$ the twisted conjugacy class $[e]_{\varphi}$ is a subgroup of $G$. Let $N$ be a normal $\varphi$-admissible subgroup of $G$. Then the class $[\overline{e}]_{\overline{\varphi}}$ is a subgroup of $G/N$.
\end{prp}

\begin{ttt}\label{cor2} Let  $G$ be a Chevalley group of type $\Phi\neq A_1$ over such a  field $F$ of zero characteristic, that the transcendence degree of $F$ over $\mathbb{Q}$ is finite. If one of the following conditions hold
 \begin{enumerate}
\item  $\Phi$ has one of types $A_l~(l\geq 7)$, $B_l~(l\geq 4)$, $E_8,~F_4,~G_2$,
\item  $\Phi$ has one of types $A_l~(l=2,3,4,5,6),~B_l~(l=2,3),~C_l~(l\geq 3),~D_l~(l\geq 4),~E_6,~E_7$ and the equality  $f(T)=a$ is solvable in $F$ for any $a$, where $f$ is a polynomial, which is defined in lemma \ref{lem20},
\end{enumerate}
then $\varphi$-conjugacy class of the unit element $[e]_{\varphi}$ is a subgroup of $G$ if and only if $\varphi$ acts identically modulo center of group $G$.
\end{ttt}
\textbf{Proof.} Sufficiency of the statement follows from the proposition \ref{pr4}. Let us prove necessity.

Let the set $[e]_{\varphi}$ be a subgroup of group $G$. Then by the proposition \ref{pr6} class $[\overline{e}]_{\overline{\varphi}}$ is a subgroup of $G/Z(G)\cong\Phi(F)$. According to proposition \ref{pr5} this subgroup is normal in $G/Z(G)\cong\Phi(F)$, and since $\Phi(F)$ is simple group, then this subgroup is trivial either coincides with $\Phi(F)$. Since (by theorem \ref{ttt1}) group $\Phi(F)$ possess $R_{\infty}$ property, then the class $[\overline{e}]_{\overline{\varphi}}$ can't coincide with group $\Phi(F)$ and therefore $[\overline{e}]_{\overline{\varphi}}=\overline{e}$. It means that $\overline{\varphi}$ is identical automorphism of group $G/Z(G)\cong\Phi(F)$, i.e. $\varphi$ is central automorphism of group $G$.\hfill $\square$

\begin{ttt}\label{cor3} Let  $G$ be a Chevalley group of type $\Phi\neq A_1$ over the field $F$ of zero characteristic. If an automorphism group of $F$ is torsion, then the twisted conjugacy class $[e]_{\varphi}$ of unit element is a subgroup of $G$ if and only if $\varphi$ is central automorphism of group $G$.
\end{ttt}
\textbf{Proof} is the same as the proof of the theorem \ref{cor2}.\hfill $\square$

\end{document}